\documentclass[11pt]{article}
\usepackage{amssymb}
\usepackage{latexsym}

\begin{document}

\title{\bf La formule de la trace pour les tissus planaires}

\author{ Jean-Paul Dufour  \\}

\date{}
\maketitle

\begin{abstract}
We give a complete proof of the fact that the trace of the curvature of the connection associated to a planar d-web ($d>3$) is the sum of the Blaschke curvature of its sub 3-webs.
\end{abstract}

\section{\bf Introduction.}

Un $d$-tissu du plan est une famille de $d$ feuilletages d'un ouvert  du plan qui sont deux \`a  deux transverses. Il y a plusieurs m\'ethodes \'equivalentes  pour donner ces feuilletages. Il y a des m\'ethodes "explicites" o\`u chaque feuilletage est donn\'e soit par les trajectoires d'un champ ou, ce qui revient au m\^eme, en se donnant les pentes $m_i(x,y)$ de chaque feuilletage dans un syst\`eme de coordonn\'ees $(x,y)$ ou bien encore en se donnant ses int\'egrales premi\`eres $f_i$ (les feuilletages sont donn\'es par les courbes de niveau des $f_i$). Alain H\'enaut a d\'evelopp\'e la m\'ethode "implicite" qui pr\'esente les feuilles comme les trajectoires d'une \'equation diff\'erentielle implicite du type $F(x,y,y')=0$ o\`u $F$ est un polyn\^ome de degr\'e $d$ en $y'$ \`a coefficients d\'ependant de $x$ et $y.$ Pour une bibliographie relativement compl\`ete sur ce domaine nous renvoyons au livre de J.V. Pereira et L. Pirio de 2015 \cite{PP} ou au texte de J. V. Periera au S\'eminaire Bourbaki  de 2007 \cite{JP}.

Dans ce travail, on consid\`ere un $d$-tissu  $W(f_1,\dots ,f_d)$ sur un ouvert $U$ du plan donn\'e (explicitement) par ses $d$ int\'egrales premi\`eres $f_1,\dots ,f_d.$ Un invariant important d'un tel tissu est son "rang", c'est  la dimension de l'espace vectoriel de ses "relations ab\'eliennes"  $\sum_{i=1}^d h_i(f_i)=0$ o\`u les $h_i$ sont des fonctions d'une variable, nulles en un point fix\'e. En 2004 et dans le contexte "implicite" Alain H\'enaut \cite{AH} a montr\'e que l'on pouvait  associer au tissu un fibr\'e vectoriel sur $U$ muni d'une connexion $\nabla$ dont l'une des propri\'et\'es est que sa courbure est nulle si et seulement si le rang du tissu a la valeur maximale possible $(d-1)(d-2)/2.$ A la m\^eme \'epoque L. Pirio a soutenu sa th\`ese \cite{P} dans laquelle, entre autre, il revisite des travaux anciens  de A. Pantazi \cite{AP}, pour construire une connexion analogue dans le cas "explicite". En 2005 Olivier Ripoll \cite{OR}, sous la direction d'A.H\'enaut, a soutenu une th\`ese sur ces sujets. Ces deux auteurs ont construit des  programmes Maple qui calculent  la connexion et sa courbure pour $d=3,$ $d=4$ et $d=5,$ dans le cas "explicite" pour L. Pirio et dans le cas "implicite" pour O.  Ripoll.  En 2007  Vincent Cavalier et Daniel Lehmann \cite{CL} ont g\'en\'eralis\'e les constructions pr\'ec\'edentes aux tissus de codimension 1 en dimension $n$ arbitraire, pourvu que ceux-ci soient "ordinaires" (ce qui est toujours le cas pour $n=2$) et tels qu'il existe un entier $k_0$ tel que $d=(n-1+k_0)!/((n-1)!k_0!)$ (ce qui est toujours le cas pour $n=2$ avec $k_0=d-1$). Travaillant dans le cas "explicite" avec les pentes des feuilletages, ils ont construit un fibr\'e vectoriel muni d'une connexion qui g\'en\'eralise celle de A.Pantazi. En 2014  D.
 Lehmann et l'auteur \cite{DL} ont reconstruit ce fibr\'e vectoriel et sa connexion \`a partir des int\'egrales premi\`eres des feuilletages et r\'edig\'e un programme Maple qui, non seulement calcule la connexion et la courbure des $d$-tissus plans pour tout $d,$ mais aussi qui peut fonctionner en toute dimension. Dans la suite on note $\nabla$ cette connexion.

 En 1933 W. Blaschke \cite{WB} avait attach\'e une "courbure" aux 3-tissus du plan. La courbure de Blaschke est  celle de $\nabla$ dans le cas particulier $d=3.$
 
Dans leurs  th\`eses et travaux suivants L. Pirio et O. Ripoll ont conjectur\'e le r\'esultat suivant.

\noindent {\bf FORMULE DE LA TRACE.} {\bf La trace de la courbure de  $\nabla$ est la somme des courbures de Blaschke des sous-3-tissus $W(f_i,f_j,f_k)$ de $W(f_1,\dots ,f_d)$.}

Nous voyons la courbure de Blaschke, et la trace de la courbure de $\nabla ,$ comme des 2-formes \`a valeurs scalaires sur le plan ; cela donne sens \`a la formule pr\'ec\'edente.

On peut voir les pr\'emisses de cette formule dans les travaux de N. Mihaileanu \cite{NM}  et A. Pantazi \cite{AP}. L. Pirio et O. Ripoll  ind\'ependamment en ont donn\'e des preuves pour $d=4,$  $5$ et $6$ Ils ont  aussi pr\'esent\'e  un plan de d\'emonstration pr\'ecis pour $d$ quelconque. L. Pirio a donn\'e \`a cette formule le nom de "formule de Mihaleanu."

Dans ce travail nous proposons une d\'emonstration compl\`ete de ce r\'esultat  bas\'ee sur la m\'ethode utilis\'ee pour construire le programme de D. Lehmann et de l'auteur \cite{DL}. Elle est un peu longue mais tous les calculs sont explicit\'es et \'el\'ementaires. La m\'ethode "explicite" a l'avantage de permettre  des raisonnements par r\'ecurrence sur le nombre $d$ car on peut comprendre  plus facilement ce qui advient de la connexion quand on rajoute une ($d+1$)-\`eme fonction. Dans le cadre "implicite" on ne travaille qu'avec les polyn\^omes sym\'etriques de ces fonctions et l'effet de l'ajout d'une nouvelle fonction est un peu plus cach\'e. Ceci dit,  le d\'etail des calculs que l'on trouvera dans certaines parties de la d\'emonstration laisse penser que l'on pourrait avoir une d\'emonstration plus simple dans le contexte "implicite".

\section{Construction de la connexion du $d$-tissu.} 

Nous rappelons, sans \'ecrire tous  les d\'etails, comment nous construisons le fibr\'e vectoriel et la connexion dans le texte de D. Lehmann et l'auteur\cite{DL}. 

On travaille au voisinage d'un point $P$ du plan et on impose que les fonctions $f_1$,...,$f_d$ soient nulles en $P$. On choisit des coordonn\'ees locales $x$ et $y$ qui s'annulent elles aussi en $P.$ 
Si $f$ est une fonction d\'efinie sur un voisinage de $P,$ $f_x$ (resp. $f_y$) d\'esigne la d\'eriv\'ee de $f$ par rapport \`a $x$ (resp. $y$). Ainsi $f_{xx}$ d\'esigne la d\'eriv\'ee seconde de $f$ par rapport \`a $x$....

Les relations ab\'eliennes de notre tissu sont les relations
$$\sum_{i=1}^{d}h_i(f_i)=0$$
o\`u les $h_i$ sont des fonctions d'une variable qui s'annulent \`a l'origine.

Pour \'etudier cette relation on  d\'erive successivement ses deux membres par rapport aux deux variables.

 On note $\omega^r_i=h_i^{(r)}(f_i)$ , o\`u $h_i^{(r)}$ d\'esigne la d\'eriv\'ee $r$-i\`eme de $h_i.$

A l'ordre 1 on a les deux \'equations 
$$\sum_{i=1}^{d}\omega^1_i{ f_x}=0$$$$\sum_{i=1}^{d}\omega^1_i{ f_y
}=0$$
que l'on r\'ecrit sous forme matricielle 
$$P_2(\omega^1_1,\dots ,\omega^1_d)=0$$
o\`u $ P_2 $ est la matrice jacobienne de $(x,y)\mapsto (f_1(x,y),\dots ,f_d(x,y)).$

A l'ordre 2 on a 3 \'equations que l'on range en prenant pour premi\`ere celle qui correspond \`a la d\'eriv\'ee $\partial\over{\partial x^2}$, la deuxi\`eme  \`a $\partial\over{\partial x\partial y}$, la troisi\`eme \`a $\partial\over{\partial y^2}$. On les \'ecrit matriciellement sous la forme
$$G^2_{3}(\omega^1_1,\dots ,\omega^1_d)+P_3(\omega^2_1,\dots ,\omega^2_d)=0$$.

Plus g\'en\'eralement, on range les  \'equations d'ordre $r-1$ en imposant l'ordre $${\partial\over{\partial x^{r-1}}},   {\partial\over{\partial x^{r-2}\partial y}},\dots , {\partial\over{\partial y^{r-1}}}$$  
pour les d\'erivations et on les
 r\'ecrit  sous la forme matricielle
$$G^{r-1}_{r}(\omega^1_1,\dots ,\omega^1_d)+\cdots +G^{2}_{r}(\omega^{r-2}_1,\dots ,\omega^{r-2}_d)+P_{r}(\omega^{r-1}_1,\dots ,\omega^{r-1}_d)=0.$$

Les matrices $G_r^j$ ont des coefficients qui sont des expressions polynomiales des d\'eriv\'ees partielles des $f_i$ d'ordre 1 \`a $j.$

Lorsque cela nous para\^itra utile pour rendre notre texte plus clair, nous rajouterons l'indice $(f_1,\dots ,f_d)$ \`a nos matrices. Ainsi on \'ecrira $G^i_{j ;f_1,\dots ,f_d}$ \`a la place de $G^i_j$ pour pr\'eciser quelles sont les fonctions en jeu.

La premi\`ere \'etape du programme donn\'e dans \cite{DL} calcule les coefficients des matrices $G_r^s$ et $P_r$ par r\'ecurrence. 
Nous retiendrons simplement que les matrices $P_r$ ont des colonnes de la forme  
$$f_i^{r-1}=\left|\matrix{( f_{ix})^{r-1} \cr ( f_{i x})^{r-2}( f_{iy})\cr .\cr .\cr ( f_{ix})^{r-j}( f_{iy})^{j-1}
\cr .\cr. \cr(f_{i y})^{r-1}\cr }\right| $$  et que, si l'on note $ G^{r}(f_i)$ les colonnes de $G_r^{2}$ et  $G^{r+}(f_i)$ la $(r-1)$-colonne obtenue en supprimant la derni\`ere composante de $G^r(f_i)$, on a 
 les relations de r\'ecurrence  
$$G^{r+}(f_i)=f_{ix}.G^{r-1}(f_i)+{\partial f^{r-2}_i\over\partial x}.$$

On construit par blocs la matrice \`a $(d+1)(d-2)/2$ lignes et $(d-2)d$ colonnes
$$MM=\left|\matrix{P_2 & 0&.&.&.&0\cr G_3^2 &P_3&0&.&.&0\cr G_4^3&G_4^2&P_4&0&.&0\cr .&.&.&.&.&.\cr .&.&.&.&.&.\cr .&.&.&.&.&.\cr G_{d-1}^{d-2}& G_{d-1}^{d-3}&.&.& G_{d-1}^{2}&P_{d-1} \cr}\right|$$
Le noyau de cette matrice donne le fibr\'e de rang $(d-1)(d-2)/2,$ et de base le plan des $(x,y),$ sur lequel la connexion $\nabla$ sera d\'efinie.

On construit $\nabla$ comme suit.

On consid\`ere d'abord  la matrice par blocs \`a $(d-2)d$ lignes et $(d-2)d$ colonnes
$$\Delta =\left|\matrix{0 & Id_d&0&.&.&0\cr 0 &0 &Id_d&0&.&0\cr .&.&.&.&.&.\cr .&.&.&.&.&.\cr .&.&.&.&.&.\cr 0&. &.&.&0&Id_d\cr A_{d-1}& A_{d-2}&.&.& A_{3}&A_{2} \cr}\right|$$
o\`u $Id_d$ d\'esigne la matrice identit\'e \`a $d$ lignes et
$$A_j=-P_{d}^{-1}.G^j_{d}.$$
 On note $D_x$ (resp. $D_y$) la matrice  diagonale dont les \'el\'ements diagonaux sont $( f_{1 x},\dots ,  
 f_{d x})$ (resp.  $( f_{1 y},\dots ,  
 f_{d y})).$ On consid\`ere maintenant la matrice $DD_x$ (resp. $DD_y$)  diagonale par blocs, \`a $(d-2)d$ lignes, dont les blocs diagonaux sont tous \'egaux \`a $D_x$ (resp. $D_y$). 

On consid\`ere la matrice carr\'ee \`a $(d-2)d$ colonnes $\Delta_x=  DD_x.\Delta$ (resp. $\Delta_y= DD_y.\Delta$).

Les coefficients des deux matrices $\Delta_x$ et $\Delta_y$  d\'ependent de $(x,y)$ ; elles donnent donc des morphismes du fibr\'e trivial de rang $(d-2)d$ sur le plan. On a une connexion $\nabla^0$ sur ce fibr\'e en prenant $\nabla^0_{\partial\over\partial x}={\partial\over\partial x}-\Delta_x $  et la m\^eme chose en rempla\c{c}ant $x$ par $y.$ On peut voir que cette connexion pr\'eserve le sous-fibr\'e donn\'e par le noyau de $MM$.

La connexion $\nabla$ est alors la restriction de $\nabla^0$ au  noyau de $MM$. 

Pour construire une base du noyau de $MM$ nous proc\'edons comme suit. 

On remarque que chacun des blocs "diagonaux" $P_r$ de $MM$ est tel que ses $r$ premi\`eres colonnes forment une sous-matrice inversible. Cela nous m\`ene \`a changer un peu l'\'ecriture de nos variables : on  r\'ecrit
 $$(\omega^1_1,...,\omega^1_d;\omega^2_1,...,\omega^2_d;...;\omega^{d-2}_1,...,\omega^{d-2}_d),$$
 plut\^ot sous la forme 
$$(\omega^1_1,\omega^1_2,\beta_3^1,...,\beta_d^1;\omega^2_1,\omega_2^2,\omega^2_3,\beta_4^2,...,\beta_d^2;...;\omega^{d-2}_1,...,\omega^{d-2}_{d-1},\beta_d^{d-2}).$$
C'est \`a dire que l'on remplace $\omega_i^r$ par $\beta_i^r$  pour $i> r+1$.

 On obtient une base 
$$B=\{e^1_3,\dots , e^1_{d};e^2_4,\dots ,e^2_{d}; \dots ;e^{d-3}_{d-1},e^{d-3}_d;e^{d-2}_d\}$$ 
du noyau de $MM$ en  prenant pour $e^r_i$ le vecteur du noyau de $MM$ dont les coordonn\'ees $\beta^s_j$  sont toutes nulles sauf $\beta^r_{i}$ qui est 1.

Notons $\Omega_x$ (resp. $\Omega_y$) les matrices de  $\nabla_{\partial\over\partial x}$ (resp. $\nabla_{\partial\over\partial y}$) par rapport \`a la base  $B$.
Ses coefficients  sont obtenus comme suit. On remarque d'abord que si l'on d\'erive n'importe quel vecteur $e^r_i$ de la base $B$ par rapport \`a $x$ ou $y$ on obtient un vecteur dont toutes les composantes $\beta^s_j$ sont nulles. Ainsi  $\nabla_{\partial\over\partial x}(e^r_i)$ a des composantes $\beta^s_j$ qui sont celles de $-\Delta_x(e^r_i).$ On note $\Omega_{x,s,j}^{ r,i}$ sa composante $\beta^s_j$ ; c'est sa composante sur le vecteur de base $e^s_j.$ Alors $\Omega_x$ est la matrice qui a les coefficients  $\Omega_{x,s,j}^{ r,i}$ ; autrement dit, on a  
 $$-\Delta_x(e^r_i)=\sum_{s=1}^{d-2}\sum_{j=s+2}^{d} \Omega_{x,s,j}^{ r,i}e^s_j.$$
 On agit de m\^eme en permutant $x$ et $y$ pour calculer $\Omega_y.$

\section{Calcul de la trace de la courbure.}

{\bf D\'efinition.} On consid\`ere le tissu $W(f_1,\dots ,f_s)$   donn\'e par ses $s$ int\'egrales premi\`eres locales $f_1,\dots ,f_s$  ($s>2$).
On lui associe les matrices $P_{s-1}=P_{s-1;f_1,\dots ,f_{s}},$ $P_{s}=P_{s;f_1,\dots ,f_{s}}$ et $G^2_s=G^{2}_{s;f_1,\dots ,f_{s}}$  d\'efinies comme dans le paragraphe pr\'ec\'edent.  La matrice $P_{s-1}$ est de rang $s-1$ et a un noyau de dimension 1 engendr\'e par un vecteur du type $(X_1,X_2,\dots , X_{s-1},1)$. On appelle {\bf \'el\'ement de trace} de $f_1,\dots ,f_s$
et on note $\gamma(f_1,\dots ,f_s)$ la derni\`ere composante du vecteur  
$$-P_{s}^{-1} G^{2}_{s}(X_1,X_2,\dots , X_{s-1},1).$$

Cet \'el\'ement de trace est aussi caract\'eris\'e par le fait que la matrice par blocs
$$M_{f_1,\dots ,f_{s}}=\left|\matrix{P_{s-1;f_1,\dots ,f_{s}} & 0 \cr G^{2}_{s;f_1,\dots ,f_{s}} & P_{s;f_1,\dots ,f_{s}}\cr }\right| $$ a un noyau engendr\'e par un vecteur de la forme
 $$(X_1,X_2,\dots , X_{s-1},1;Y_1,Y_2,\dots , Y_{s-1},\gamma(f_1,\dots ,f_s)).$$

Etudions maintenant l'expression de la courbure de $\nabla$ dans la base donn\'ee dans la section pr\'ec\'edente. Elle a la matrice
$$K={\partial\over\partial y}(\Omega_x)-{\partial\over\partial x}(\Omega_y)+\Omega_x.\Omega_y -\Omega_y.\Omega_x.$$
Le commutateur $\Omega_x.\Omega_y -\Omega_y.\Omega_x$ a une trace nulle, donc la trace de la courbure est la trace de la matrice 
$$KK={\partial\over\partial y}(\Omega_x)-{\partial\over\partial x}(\Omega_y).$$

\noindent{\bf Remarque.} {\sl Dans un prochain travail avec D. Lehmann  nous prouverons que la matrice $K$ a toutes ses lignes nulles sauf la derni\`ere. Donc  sa trace se r\'eduit au seul terme diagonal sur la cette derni\`ere ligne. On aurait pu penser que ce r\'esultat serait un ingr\'edient essentiel de toute preuve de la formule de la trace. Dans la d\'emonstration suivante nous proc\'edons autrement, en n'utilisant que $KK.$}

 Pour calculer la trace de $K$, il suffit donc de calculer la trace de $\Omega_x,$ de la d\'eriver par rapport \`a $y$, puis de retrancher ce que l'on obtient en \'echangeant les r\^oles de $x$ et $y$.
  
  La trace de $\Omega_x$ est la somme des $\Omega^{r,i}_{x, r,i}.$ 

Dans la suite de cette section on calcule   $\Omega^{r,i}_{x, r,i}$  pour $r$ et $i $ fix\'es. Pour cela on rappelle que $e^r_i$ est le $(d-2)d$-vecteur
$$(\omega^1_1,\omega^1_2,\beta_3^1,...,\beta_d^1;\omega^2_1,\omega_2^2,\omega^2_3,\beta_4^2,...,\beta_d^2;...;\omega^{d-2}_1,...,\omega^{d-2}_{d-1},\beta_d^{d-2})$$ du noyau de $MM$ dont les coordonn\'ees $\beta^s_j$  sont toutes nulles sauf $\beta^r_{i}$ qui est 1. La forme "triangulaire inf\'erieure par blocs" de $MM$ implique que les $\omega^s_i$  sont tous nuls pour $s<r.$ C'est dire que $e^r_i$ a la forme
$$(0,...,0;...;0,...,0;\omega^r_1,...,\omega^r_{r+1},0,...,0,1,0,...,0;\omega^{r+1}_1,...,\omega^{r+1}_{r+2},0,...,0;\omega^{r+2}_1...), $$
o\`u le 1 est \`a la $i$-\`eme place entre les deux points virgules qui l'encadrent.

On va distinguer deux cas :

\noindent {\bf Le cas $r<d-2$.}
Etudions $-\Delta_x(e_i^r)$ ; c'est, avec les notations de la section pr\'ec\'edente, $- DD_x.\Delta(e^r_j)$. La structure des lignes par blocs $Id_d$ de la partie  sup\'erieure de $\Delta$ fait que l'on a
$$-\Delta_x(e_i^r)=(\theta^1_1,...,\theta^1_d;....;\theta^ {d-2}_1,...,\theta^ {d-2}_d)$$
o\`u les $\theta^s_i$ sont tous nuls pour $s<r-1$,
 $$(\theta^{r-1}_1,...,\theta^{r-1}_d)=(-f_{1x}\omega^r_1,...,-f_{(r-1)x}\omega^r_{r+1},0,...,0,-f_{ix},0,...,0)$$ (ces deux conditions n'ayant de sens  que pour $r>1$) et
$$(\theta^{r}_1,...,\theta^{r}_d)=(-f_{1x}\omega^{r+1}_1,...,-f_{(r+2)x}\omega^{r+1}_{r+2},0,...,0).$$
On en d\'eduit
$$\Omega^{r,i}_{x, r,i}=0$$
 si $i>r+2$ et
$$\Omega^{r,r+2}_{x, r,r+2}=-f_{(r+2)x}\omega^{r+1}_{r+2}.$$
Exprimons maintenant ce qu'est $\omega^{r+1}_{r+2}$

Revenant un peu en arri\`ere (avec $i=r+2$) nous remarquons que le fait que $e_{r+2}^r$ soit dans le noyau de $MM$ et que ses composantes $\beta^s_j$ soient nulles pour $j>r+2$ et $s=r$ ou $s=r+1$ nous donne  la relation 
$$M_{f_1,...,f_{r+2}}(\omega^r_1,...,\omega^r_{r+1},1;\omega^{r+1}_1,...,\omega^{r+1}_{r+2})=0.$$
D'apr\`es la d\'efinition de d\'ebut de cette section cela veut dire que
$\omega^{r+1}_{r+2}$ {\bf est l'\'el\'ement de trace} $\gamma(f_1,\dots ,f_{r+2})$.
0n en d\'eduit finalement

$$\Omega^{r,r+2}_{x, r,r+2}=-f_{(r+2)x}\gamma(f_1,\dots ,f_{r+2}).$$

\noindent  {\bf Le cas $r=d-2$.}

 On a
$e^{d-2}_d=(0;...;0;X)$ avec $X=(\omega^{d-2}_1,\cdots ,\omega^{d-2}_{d-1},1)$ et $P_{d-1}(X)=0 $  (notation de la section pr\'ec\'edente). On a
  $$-\Delta_x (e^{d-2}_d)=(0;...;0;-D_x(X);-D_x.A_{2}(X).$$
Or, par d\'efinition, la derni\`ere  composante de $A_2(X)=-P^{-1}G_2(X)$ est l'\'el\'ement de trace  $\gamma(f_1,\dots ,f_{d})$. On en tire une formule analogue \`a celle des cas pr\'ec\'edents 

$$\Omega^{d-2,d}_{x, d-2,d}=-f_{dx}\gamma(f_1,\dots ,f_{d}).$$

On a le r\'esultat analogue lorsque l'on \'echange $x$ et $y$ et on en tire facilement
  la proposition suivante.

\noindent {\bf Proposition.} {\sl La trace de la courbure du tissu} $W(f_1,\dots , f_d)$ {est la 2-forme}
$$-\sum^d_{r=3} df_{r}\wedge d\gamma (f_1,\dots ,f_r).$$

\section{Le cas des 3-tissus.}

On consid\`ere d'abord  un 3-tissu $W(f,g,h)$ du plan. Trivialement sa courbure classique de Blaschke est aussi la trace de la courbure de la connexion associ\'ee. Suivant la proc\'edure d\'ecrite dans la section pr\'ec\'edente, on l'obtient en calculant d'abord  l'\'el\'ement de courbure $\gamma (f,g,h)$. Nous allons donner son expression pr\'ecise dans le cas particulier o\`u $h(x,y)=y$ et $f_x$ et $g_x$ ne s'annulent pas.

On adopte les notations $m_f=f_y/f_x$ et  $m_g=g_y/g_x$. On prend
$$X_1={1\over {(m_g-m_f)f_x}},\ \ X_2={1\over {(m_f-m_g)g_x}}$$
 et on voit que $(X_1,X_2,1)$ engendre le noyau de
$$P_{2;f,g,h}=\left | \matrix{f_{x} & g_{x} & 0  \cr f_{y} & g_{y} & 1 \cr}\right | .$$

Alors $\gamma (f,g,h)$ est la derni\`ere composante du vecteur
$$- P_{3;f,g,h}^{-1}G^{2}_{3;f,g,h}(X_1,X_2,1)$$

en prenant :

$$ P_{3;f,g,h}=\left | \matrix{(f_{x})^2 & (g_{x})^2 & 0  \cr f_x f_{y} & g_x g_{y} & 0 \cr (f_{y})^2 & (g_{y})^2 & 1  \cr}\right | ,$$
$$G^{2}_{3;f,g,h}=\left | \matrix{f_{xx} & g_{xx} & 0  \cr f_{xy} & g_{xy} & 0 \cr f_{yy} & g_{yy} & 0 \cr }\right |.$$

Or le m\^eme argument que celui qui permet de calculer la matrice inverse d'une matrice de Vandermonde, montre que la derni\`ere ligne de $ P_{3;f,g,h}^{-1}$ est  $(m_fm_g,-(m_f+m_g),1)$ ; on en tire
$$\gamma (f,g,h)=-(m_fm_g,-(m_f+m_g),1).G^{2}_{3;f,g,h}(X_1,X_2,1),$$
ce qui m\`ene facilement \`a

$$\gamma (f,g,h)={1\over {m_f-m_g}}\{ m_fm_g(f_{xx}/f_x-g_{xx}/g_x)-$$$$(m_f+m_g)(f_{xy}/f_x-g_{xy}/g_x)+(f_{yy}/f_x-g_{yy}/g_x)\}.$$

On en tire une formule explicite pour la courburede Blaschke de notre tissu :$$ -dy\wedge d\gamma (f,g,y).$$

\section{La m\'ethode de d\'emonstration de la formule  de la trace.}

 Pour le tissu $W(f_1,\dots ,f_s)$ la somme des courbures des sous 3-tissus est 
$$SC(f_1,\dots ,f_s)=-\sum_{0<i<j<r\leq s}df_r\wedge d\gamma(f_i,f_j ,f_r).$$ 
On note
$$Tr(f_1,\dots ,f_s)(=-\sum^s_{r=3}( df_{r}\wedge d\gamma (f_1,\dots ,f_r))$$
la trace de la courbure de la connexion $\nabla$ associ\'ee comme plus haut.
La formule de la trace dit que l'on a $$SC(f_1,\dots ,f_s)=Tr(f_1,\dots ,f_s)$$
pour tout $s$ plus grand que 3.

La formule est triviale pour $s=3.$ Nous la d\'emontrons par r\'ecurrence sur $s$ ; pour cela nous la supposons montr\'ee \`a l'ordre $d-1$ et nous allons la prouver \`a l'ordre $d.$ Comme $SC(f_1,\dots ,f_d)$ et $Tr(f_1,\dots ,f_d)$ 
sont des quantit\'es qui ne d\'ependent pas des coordonn\'ees  on peut  choisir ces coordonn\'ees pour avoir $f_d =y$ et la non-nullit\'e des $f_{ix}$ pour $i$ variant de 1 \`a $d-1.$

La formule \`a l'ordre $d-1$ nous donne
$$\sum_{r=3}^{d-1} (df_{r}\wedge d(\sum_{0<i<j<r}\gamma(f_i,f_j ,f_r))=\sum^{d-1}_{r=3}( df_{r}\wedge d\gamma (f_1,\dots ,f_r)).$$

Pour prouver la formule \`a l'ordre $d$ il suffit de prouver la relation
$$\gamma(f_1,\dots ,f_d)=\sum_{0<i<j<d}\gamma(f_i,f_j ,f_d).$$
 C'est ce que nous allons faire dans la suite de ce travail en utilisant l'hypoth\`ese simplificatrice $f_d=y.$

Les quantit\'es $\gamma(f_1,\dots ,f_d)$ et $\sum_{0<i<j<d}\gamma(f_i,f_j ,f_d)$ ont des expressions lin\'eaires dans les d\'eriv\'ees secondes $f_{sxx},$ $f_{sxy}$ et $f_{syy}$ avec des coefficients qui ne d\'ependent que des d\'eriv\'ees premi\`eres.
On va voir que ces coefficients sont les m\^emes dans les deux expressions.

\section{La somme des courbures des sous 3-tissus.} 

La section 4 montre que l'on a la formule
$$\sum_{0<i<j<d}\gamma(f_i,f_j ,y)=\sum_{0<i<j<d}{1\over {m_i-m_j}}\{ m_im_j(f_{ixx}/f_{ix}-f_{jxx}/f_{jx})-$$$$(m_i+m_j)(f_{ixy}/f_{ix}-f_{jxy}/f_{jx})+(f_{iyy}/f_{ix}-f_{jyy}/f_{jx})\}$$
avec la notation $m_k=f_{ky}/f_{kx}.$
On en d\'eduit un d\'eveloppement
 $$\sum_{0<i<j<d}\gamma(f_i,f_j ,y)=\sum_{s<d}A_sf_{sxx}+B_sf_{sxy}+C_sf_{sxx}$$
avec
$$A_s={1\over f_{sx}}\sum_{j\neq s}{m_sm_j\over{m_s-m_j}}$$
$$B_s={-1\over f_{sx}}\sum_{j\neq s}{m_s+m_j\over{m_s-m_j}}$$
$$C_s={1\over f_{sx}}\sum_{j\neq s}{1\over{m_s-m_j}}.$$
Dans les sections suivantes nous allons montrer que $\gamma(f_1,\dots ,f_d)$ a le m\^eme d\'eveloppement.

\section{Calcul du noyau de $P_{d-1}.$}

Il nous faut calculer l'\'el\'ement de trace $\gamma(f_1,\dots ,f_d)$ (avec $f_d=y$). Si l'on revient sur sa d\'efinition, donn\'ee dans la section 3, il nous faut d'abord calculer le vecteur $(X_1,\dots ,X_{d-1},1)$ qui engendre le noyau de $P_{d-1}$.

En tenant compte du fait que $f_d=y,$ le syt\`eme $P_{d-1} (X_1,\dots ,X_{d-1},1)=0$ se r\'ecrit sous la forme
$$\sum_{i=1}^{d-1}f_{ix}^{d-1-j}f_{iy}^{j-1}X_i=0$$
$$\sum_{i=1}^{d-1}f_{iy}^{d-2}X_i=-1$$
o\`u $j$ varie de 1 \`a $d-1.$ Si l'on pose $Y_i=f_{ix}^{d-2}X_i$ et $m_i=f_{iy}/f_{ix}$ les $d-1$ premi\`eres \'equations donnent un syst\`eme matriciel $VM((Y_1,\dots ,Y_{d-1}))=(0,\dots,0,-1)$ o\`u $VM$ est une matrice de Vandermonde dont le coefficient sur la $s$-i\`eme ligne et la $r$-i\`eme colonne est $m_r^{s-1}.$ On en tire que $(Y_1, \dots ,Y_{d-1})$ est l'oppos\'e de la transpos\'ee de la derni\`ere colonne de $VM^{-1}.$

On a alors
$$Y_i={-1\over \prod_{j\neq i}(m_i-m_j)}$$
et donc
$$X_i={-1\over{f_{ix}^{d-2} \prod_{j\neq i}(m_i-m_j)}}$$
pour $i$ variant de 1 \`a $d-1.$

\section{Calcul de la matrice $G^2_d$.}

Dans la section 2 nous avions not\'e $G^d(f_1),\dots ,G^d(f_{d}),$ ses colonnes (avec $f_d=y$). 

On remarque d'abord que, pour une fonction arbitraire $f$, le $i$-\`eme coefficient de  $G^d(f)$ est de la forme
$$G^d_i(f)=a_i^df_x^{d-i-2}f_y^{i-1}f_{xx}+b_i^df_x^{d-i-1}f_y^{i-2}f_{xy}+c_i^df_x^{d-i}f_y^{i-3}f_{yy},$$
avec la convention d'\'ecriture  que les puissances n\'egatives des $f_x$ ou $f_y$ sont nulles ; les $a_i^d,$ $b_i^d$ et $c_i^d$ sont des nombres que nous allons d\'eterminer.

Pour des raisons de sym\'etrie par rapport aux deux d\'erivations $\partial /\partial x$ et  $\partial /\partial y$, on a les relations suivantes :
$$a_i^d=c^d_{d-i+1},\ \ \ b_i^d=b^d_{d-i+1}.$$

Nous utilisons la relation de r\'ecurrence
$$G^{d+}(f)=f_x.G^{d-1}(f)+{\partial f^{d-2}_i\over\partial x}$$
que nous avions donn\'ee en section 2. Elle nous donne
$$G^d_i(f)=f_xG^{d-1}_i(f)+(d-1-i)f_x^{d-i-2}f_y^{i-1}f_{xx}+(i-1)f_x^{d-i-1}f_y^{i-2}f_{xy}$$
pour $i<d.$

On en tire les relations de r\'ecurrence :
$$a_i^d=a_i^{d-1}+d-1-i$$
$$b_i^d=b_i^{d-1}+i-1$$
$$c_i^d=c_i^{k-1}$$
pour $i<d.$

On a aussi les relations \'evidentes :
$$a_1^3=1,\ \ b_1^3=c_1^3=0$$
$$a_2^3=0,\ \ b_2^3=1,\ \ c_1^3=0$$
$$a_1^3= b_1^3=0,\ \ c_1^3=1.$$

Utilisant ces relations, les relations de sym\'etrie et de r\'ecurrence ci-dessus on obtient :

$$a_i^d={(d-1-i)(d-i)\over 2}$$
$$b_i^d=(i-1)(d-i)$$
$$c_i^d={(i-2)(i-1))\over 2}$$
pour $ i $ compris entre 1 et $d.$

\section{La derni\`ere ligne de $P_d^{-1}$.}

Pour calculer l'\'el\'ement de trace $\gamma(f_1,\dots , f_{d-1},y)$ nous aurons besoin d'un autre ingr\'edient :
la derni\`ere ligne $\alpha =(\alpha_1,\dots , \alpha_d)$ de $P_d^{-1}$. Cette ligne est caract\'eris\'ee par le fait que le produit de cette ligne avec chacune des $r-1$ premi\`eres colonnes de $P_d$ est nul et son produit avec la derni\`ere colonne est 1. On a donc les \'equations
$$\sum_{j=1}^{d-1}\alpha_jf_{ix}^{d-j}f_{iy}^{j-1}=0$$
pour tout $i$ variant de 1 \`a $d-1$ et
$$\alpha_d=1.$$ 
En divisant les deux membres des $d-1$ premi\`eres \'equations par $f_{ix}$ on peut les remplacer par
$$\sum_{j=1}^{d-1}\alpha_jm_i^{j-1}=0.$$
Comme lorsque l'on calcule l'inverse d'une matrice de Vandermonde, on introduit le polyn\^ome $P(t)=\sum_{j=1}^{d-1}\alpha_j t^{j-1}$ et les relations pr\'ec\'edentes montrent que ce polyn\^ome admet les racines $m_1,\dots ,m_{ d-1}$ et 1 comme coefficient du terme de plus haut degr\'e.
On en d\'eduit
$$\alpha =((-1)^{d-1}S_{d-1},(-1)^{d-2}S_{d-2}, \dots ,-S_1,1),$$
o\`u les $S_i$ sont les polyn\^omes sym\'etriques en $m_ 1,m_2,\dots ,m_{d-1}.$

\section{Le calcul de l'\'el\'ement de trace $\gamma (f_1,\dots ,f_{r-1},y)$.}

Rappelons que, par d\'efinition, $\gamma (f_1,\dots ,f_{r-1},y)$ est le dernier coefficient de $-P_k^{-1}G^2_d(X_1,\dots ,X_{d-1},1)$ o\`u les $X_i$ sont ceux de la section 7 ;  donc c'est
 le produit scalaire usuel des deux vecteurs $\alpha$ (voir section 9) et $G^2_d(X_1,\dots ,X_{d-1},1)$. On \'ecrit ce r\'esultat sous la forme
$$\gamma (f_1,\dots ,f_{r-1},y)=-\alpha . G^2_d(X_1,\dots ,X_{d-1},1).$$ 

On en d\'eduit 
$$\gamma (f_1,\dots ,f_{r-1},y)=-\sum_{s=1}^{d-1}X_s(\alpha .G^d(f_s))$$
en remarquant que $G^d(y)$ est nulle.

On rappelle la formule de la section 8 :
$$G^d_i(f_s)=a_i^df_{sx}^{d-i-2}f_{sy}^{i-1}f_{sxx}+b_i^df_{sx}^{d-i-1}f_{sy}^{i-2}f_{sxy}+c_i^df_{sx}^{d-i}f_{sy}^{i-3}f_{syy}$$
pour la $i$-\`eme composante de $G^d(f_s).$ On rappelle que, pour cette formule, les puissances n\'egatives de d\'eriv\'ees de fonctions sont nulles par convention. On a donc la formule
$$\alpha .G^d(f_s)=a^sf_{sxx}+b^sf_{sxy}+c^sf_{syy},$$
avec
$$a^s=\sum_{i=1}^d\alpha_ia_i^df_{sx}^{d-i-2}f_{sy}^{i-1},$$
$$b^s=\sum_{i=1}^d\alpha_ib_i^df_{sx}^{d-i-1}f_{sy}^{i-2},$$
$$c^s=\sum_{i=1}^d\alpha_ic_i^df_{sx}^{d-i}f_{sy}^{i-3}.$$

Plus pr\'ecis\'ement, on a donc
$$a^s=\sum_{i=1}^d (-1)^{d-i}{{(d-i-1)(d-i)}\over 2}S_{d-i}f_{sx}^{d-i-2}f_{sy}^{i-1},$$
$$b^s=\sum_{i=1}^d (-1)^{d-i}(i-1)(d-i)S_{d-i}f_{sx}^{d-i-1}f_{sy}^{i-2},$$
$$c^s=\sum_{i=1}^d (-1)^{d-i}{{(i-1)(i-2)}\over 2}S_{d-i}f_{sx}^{d-i}f_{sy}^{i-3}$$
ou encore
$$a^s=f_{sx}^{d-3}\sum_{i=1}^d (-1)^{d-i}{{(d-i-1)(d-i)}\over 2}S_{d-i}m_s^{i-1},$$
$$b^s=f_{sx}^{d-3}\sum_{i=1}^d (-1)^{d-i}(i-1)(d-i)S_{d-i}m_s^{i-2},$$
$$c^s=f_{sx}^{d-3}\sum_{i=1}^d (-1)^{d-i}{{(i-1)(i-2)}\over 2}S_{d-i}m_s^{i-3}.$$

On a la formule
$$S_{d-i} =m_sS^s_{d-i-1}+S^s_{d-i},$$
en notant $S^s_k$ le $k$-i\`eme polyn\^ome sym\'etrique dans les variables $$m_1, \dots ,m_{s-1},m_{s+1},\dots ,m_{d-1}$$ avec, par convention, $S^s_j=0$  pour $j<0$ ou $j>d-2.$
C'est dire l'on oublie $m_s$ dans les $S^s_j.$ En portant cela dans \'equations pr\'ec\'edentes on obtient
$$a^s=f_{sx}^{d-3}\sum_{i=1}^d (-1)^{d-i}{{(d-i-1)(d-i)}\over 2}(S_{d-i}^sm_s^{i-1}+S^s_{d-i-1}m_s^i),$$
$$b^s=f_{sx}^{d-3}\sum_{i=1}^d (-1)^{d-i}(i-1)(d-i)(S_{d-i}^sm_s^{i-2}+S^s_{d-i-1}m_s^{i-1}),$$
$$c^s=f_{sx}^{d-3}\sum_{i=1}^d (-1)^{d-i}{{(i-1)(i-2)}\over 2}(S^s_{d-i}m_s^{i-3}+S^s_{d-i-1}m_s^{i-2}).$$
En r\'eordonnant les termes en fonction des puissances de $m_s,$ on arrive \`a
$$a^s=f_{sx}^{d-3}\sum_{i=1}^d (-1)^{d-i}(d-i-1)S_{d-i-1}^sm_s^i,$$
$$b^s=f_{sx}^{d-3}\sum_{i=1}^d (-1)^{d-i}(d-2i-2)S_{d-i-2}^sm_s^i,$$
$$c^s=f_{sx}^{d-3}\sum_{i=1}^d (-1)^{d-i-1}(i+1)S^s_{d-i-3}m_s^{i}.$$

\section{Fin de la preuve de la formule de la trace.}

On a la relation
 $$\gamma (f_1,\dots ,f_{d-1},y)=-\sum_{s=1}^{d-1}X_s(a^sf_{sxx}+b^sf_{sxy}+c^sf_{syy}).$$
Donc nous aurons montr\'e la formule de la trace si l'on prouve les relations
$$X_sa^s=-A_s,\ \ X_sb^s=-B_s,\ \ X_sc^s=-C_s,$$
les $A_s$, $B_s$ et $C_s$ \'etant ceux d\'efinis dans la section 6. Il revient au m\^eme de prouver les relations
$$a^s=-A_s/X_s,\ \  b^s=-B_s/X_s,\ \ c^s=-C_s/X_s.$$
Pour simplifier les notations on ne d\'emontrera ces relations que dans le cas $s=d-1$ car exactement la m\^eme m\'ethode fonctionne dans les autres cas. Toujours pour simplifier, nous \'ecrirons $m$ \`a la place de $m_{d-1}.$ 
On a 
$$-A_{d-1}/X_{d-1}={1\over{ f_{(d-1)x}^{d-3}}}(\prod_{r=1}^{d-2}(m-m_r))\sum_{j=1}^{d-2}{{mm_j}\over{m-m_j}}.$$
On a la relation
$$(\prod_{r=1}^{d-2}(m-m_r))\sum_{j=1}^{d-2}{{mm_j}\over{m-m_j}}=$$$$\sum_{j=1}^{d-2}(m-m_1)\cdots (m-m_{j-1})mm_j(m-m_{j+1})\cdots (m-m_{d-2})$$
Pour calculer cette quantit\'e on introduit la  fonction $$P(t)=\prod_{j=1}^{d-2}(t-m_j)$$
et il est facile de voir que l'on la relation
$$(\prod_{r=1}^{d-2}(t-m_r))\sum_{j=1}^{d-2}{{tm_j}\over{t-m_j}}=t(tP'(t)-(d-2)P(t)).$$
Mais $P(t)$ est aussi le polyn\^ome
$$P(t)=\sum_{j=1}^{d-2}(-1)^{d-2-j}S^{d-1}_{d-2-j}t^j,$$
o\`u  $S^{d-1}_r$ est le polyn\^ome sym\'etrique de degr\'e $r$ dans les variables $m_1,$ ...,$m_{d-2}.$
Cela m\`ene \`a 
$$(\prod_{r=1}^{d-2}(t-m_r))\sum_{j=1}^{d-2}{{tm_j}\over{t-m_j}}=\sum_{i=1}^{d-2}(-1)^{d-i}(d-i-1)S_{d-i-1}^{d-1}t^i,$$
et, en posant $t=m$ et rajoutant le facteur ${1\over{ f_{(d-1)x}^{d-3}}},$ on obtient
$$a^{d-1}=-A_{d-1}/X_{d-1}, $$
et, par la m\^eme m\'ethode, \`a 
$$a^{s}=-A_{s}/X_{s}, $$
pour tout $s.$ 

Par ailleurs, on a 
$$-B_{d-1}/X_{d-1}={1\over{ f_{(d-1)x}^{d-3}}}(\prod_{r=1}^{d-2}(m-m_r))\sum_{j=1}^{d-2}{{m+m_j}\over{m-m_j}}.$$
Comme pour le calcul pr\'ec\'edent on obtient facilement
$$(\prod_{r=1}^{d-2}(t-m_r))\sum_{j=1}^{d-2}{{t+m_j}\over{t-m_j}}=2tP'(t)-(d-2)P(t)$$
qui, en revenant \`l'expression polynomiale de $P(t),$ m\`ene \`a
$$(\prod_{r=1}^{d-2}(t-m_r))\sum_{j=1}^{d-2}{{t+m_j}\over{t-m_j}}=\sum_{i=1}^{d-2}(-1)^{d-i-1}(d-2i+2)S_{d-i-2}^{d-1}t^i,$$
d'o\`u il d\'ecoule
$$b^{d-1}=-B_{d-1}/X_{d-1}$$
et de la m\^eme mani\`ere
$$b^{s}=-B_{s}/X_{s}$$
pour tout $s.$

Enfin, on a

$$-C_{d-1}/X_{d-1}={1\over{ f_{(d-1)x}^{d-3}}}(\prod_{r=1}^{d-2}(m-m_r))\sum_{j=1}^{d-2}{1\over{m-m_j}}.$$
Comme 
$$(\prod_{r=1}^{d-2}(t-m_r))\sum_{j=1}^{d-2}{1\over{t-m_j}}=P'(t),$$
on a
$$(\prod_{r=1}^{d-2}(t-m_r))\sum_{j=1}^{d-2}{1\over{t-m_j}}=\sum_{j=1}^{d-2}(-1)^{d-2-j}jS^{d-1}_{d-2-j}t^{j-1},$$
et
$$c^{d-1}=-C_{d-1}/X_{d-1}$$
puis, par la m\^eme m\'ethode,
$$c^{s}=-C_{s}/X_{s},$$
pour tout $s.$

Ceci ach\`eve notre d\'emonstration de la formule de la trace pour les tissus planaires.

\section{Exemple d'application de la formule de la trace.}

Alain H\'enaut a propos\'e la conjecture suivante. {\sl Soit $W$ un $d$-tissu du plan donn\'e de mani\`ere implicite par le polyn\^ome
$$F=(y')^{d}+f(x,y)\ ;$$
alors la courbure de la connexion associ\'ee est nulle si et seulement si $f(x,y)$ est d\'ecomposable, c'est \`a dire de la forme $X(x)Y(y).$} 

Remarquons que dans l'ouvert o\`u $F$ est non nul, les $d$ racines de $F$ sont toutes de la forme $m_i=R(x,y)\lambda_i$  o\`u les $\lambda_i$ sont des constantes deux \`a deux diff\'erentes et la fonction $R(x,y)$ est ind\'ependante de l'indice  $i.$ De plus $f(x,y)$ est d\'ecomposable si et seulement si $R(x,y)$ l'est. 

 Alors la conjecture d'H\'enaut est un corollaire du lemme suivant.

{\bf Lemme.} {\sl Soit $W$ un $d$-tissu plan donn\'e par ses pentes $m_i$ pour $i$ variant de 1 \`a $d$. On suppose que l'on a 
$$ m_i=R(x,y)\lambda_i$$
 o\`u les $\lambda_i$ sont des constantes deux \`a deux diff\'erentes et la fonction $R(x,y)$ est ind\'ependante de l'indice  $i$ et non nulle. 
Alors la courbure de la connexion associ\'ee est nulle si et seulement si $R(x,y)$ est d\'ecomposable.}

Nous allons donner une preuve de ce lemme en montrant d'abord  le sens direct : si $R(x,y)$ est de la forme $X(x)Y(y)$ alors la courbure est nulle puis la r\'eciproque.

{\bf 1- On suppose} $R(x,y)=A(x)B(y).$

Les feuilles des $d$ feuilletages sont les trajectoires des $d$ champs de vecteurs
$$X_i={\partial\over\partial x}+A(x)B(y)\lambda_i{\partial\over\partial y}.$$
Ce sont aussi les trajectoires de
$$Y_i=1/A(x){\partial\over\partial x}+\lambda_iB(y){\partial\over\partial y}.$$
Or des changements de la coordonn\'ee $x,$ d'une part, et $y,$ d'autre part, permettent de rectifier les champs de vecteurs $1/A(x){\partial\over\partial x}$ et $B(y){\partial\over\partial y}.$ Ces changements nous ram\`enent au cas o\`u tous les $Y_i$ sont \`a coefficients constants, donc au cas o\`u les $d$ feuilletages sont tous form\'es de segments parall\`eles. Or on sait que ces tissus
sont \`a courbure nulle.

{\bf 2- La r\'eciproque.}
On suppose que  $W$ est de courbure nulle. Alors la trace de cette courbure est encore nulle, donc, par la {\bf formule de la trace}, la somme des courbure des sous 3-tissus de $W$ est nulle.
Etudions le sous 3-tissu de $W$ donn\'e par les pentes $m_i,$ $m_j$ et $m_k.$ Un calcul \'el\'ementaire (qui peut \^etre fait par Maple) montre que sa coubure est 
 $${\partial^2\over {\partial x\partial y}}ln(R(x,y)).$$
Ainsi on voit que la somme des courbures des sous 3-tissus est nulle si et seulement si ${\partial^2\over {\partial x\partial y}}ln(R(x,y))=0,$ ou, ce qui est \'equivalent, que la fonction $R(x,y)$ est d\'ecomposable. Ceci ach\`eve la d\'emonstration.


\begin{thebibliography}{dango 9999}

\bibitem[WB]{WB} W. Blaschke, {\it Uber die Tangenten einer ebenen Kurve funfter Klasse. } Abh. Math. Semin. Hamb. Univ. 9 (1933) 313-317.

\bibitem[CL]{CL} V. Cavalier, D. Lehmann, {\it Ordinary holomorphic webs of codimension one. } arXiv 0703596v2 [mathsDS], 2007, et Ann. Sc. Norm. Super. Pisa, cl. Sci (5), vol XI (2012), 197-214.
   .
\bibitem[DL]{DL} J. P. Dufour, D. Lehmann, {\it Calcul explicite de la courbure des tissus calibr\'es ordinaires } arXiv 1408.3909v1 [mathsDG],18/08/2014.

\bibitem[AH]{AH} A. H\'enaut, {\it Planar web geometry through abelian relations  and connections}
   Annals of Math. 159 (2004)  425-445.

\bibitem[NM]{NM} N. Mihaileanu. {\it Sur les tissus plans de premi\`ere esp\`ece.} Bull. Math. Soc. Roum.Sci. 43 (1941), 23-26.
  
\bibitem[P]{P} L. Pirio, {\it Equations Fonctionnelies Ab\'eliennes et G\'eom\'etrie des tissus}
  Th\`ese de doctorat de l'Universit\'e Paris VI, 2004.

\bibitem[AP]{AP} A. Pantazi. {\it Sur la d\'etermination du rang d'un tissu plan.} C.R. Acad. Sc. Roumanie 4 (1940), 108-111.

\bibitem[PP]{PP} J.V. Pereira, L. Pirio,{\it An Invitation to Web Geometry}
  Series IMPA Monographs Vol.2, Springer  (2015).

\bibitem[JP]{JP} J.V. Pereira, {\it Algebraization of codimension one webs}
   S\'eminaire Bourbaki, 59\`eme ann\'ee, 2006-2007, $n^0$974  (mars 2007).


\bibitem[OR]{OR} O. Ripoll, {\it G\'eom\'etrie des tissus du plan et \'equations diff\'erentielles}
  Th\`ese de doctorat de l'Universit\'e de Bordeaux 1, 2005.

\end{thebibliography}
\end{document}